\documentclass{elsart}

\usepackage{amssymb,amsmath,mathrsfs,eufrak,mathabx}

\newtheorem{exm}{Example}[section]

\newcommand{\ar}{\rightarrow}
\newcommand{\la}[1]{\it {#1}\rm}

\newcommand{\mT}{\mathcal{T}}
\newcommand{\mL}{\mathcal{L}}

\newcommand{\mA}{\mathcal{A}}

\newcommand{\mC}{\mathcal{C}}

\newcommand{\mF}{\mathcal{F}}
\newcommand{\mM}{\mathcal{M}}

\newcommand{\lm}[1]{\begin{lem} #1 \end{lem}}
\newcommand{\te}[1]{\begin{thm} #1 \end{thm}}

\newcommand{\ps}[1]{\begin{prop} #1 \end{prop}}

\newcommand{\FF}{{\mathsf F }}

\begin{document}

\begin{frontmatter}



\title{Universal Algebra and Mathematical Logic}


{\large Zhaohua Luo}
\author{}


\ead{zluo@azd.com}
\ead[url]{http://www.algebraic.net/cag}



\end{frontmatter}


In this paper, first-order logic is interpreted in the framework of universal algebra, using the clone theory developed in three previous papers \cite{luo:1} \cite{luo:2} and \cite{luo:3}. 

We first define the free clone $\mT(\mL, \mC)$ of terms of a first order language $\mL$ over a set $\mC$ of parameters in a standard way. The free right algebra $\mF(\mL, \mC)$ of formulas over the clone $\mT(\mL, \mC)$ of terms is then generated by atomic formulas vis a binary operation $\Rightarrow$ and a unary operation $\forall $. The classical unary operations $\forall x_1, \forall x_2, ...$ are derived from $\forall$ via substitutions.  Structures for $\mL$ over $\mC$ are represented as perfect valuations of $\mF(\mL, \mC)$, and theories of $\mL$ are represented as filters of $\mF(\mL, \emptyset)$. Finally Godel's completeness theorem and first incompleteness theorem  are stated as expected.

A \la{(first order) language} is a nonempty set $\mL$ consisting of \la{$n$-ary function symbols} and \la{$n$-ary predicate symbols} for each $n \ge 0$. We assume $\mL$ contains a $0$-ary predicate symbol $\FF$. We say $\mL$ is a \la{language with equality} if $\mL$ contains a  $2$-ary predicate symbol $\approx$.

Let $X = \{x_1, x_2, ...\}$ be a fixed set of \la{variables}.  Let $\mC$ be a (possibly empty) set of \la{parameters}. The terms $\mT(\mL, \mC)$ of $\mL$ over $\mC$ form the smallest set of expressions containing $0$-ary function symbols, variables and parameters, which is closed under the formation rule: if $t_1, ..., t_n$ are terms of $\mL$ and if $f \in \mL$ is an $n$-ary function symbol, then the  expression $f(t_1, ..., t_n)$ is a term of $\mL$.

An \la{atomic formula} of $\mL$ over $\mC$ is either a $0$-ary predicate symbol, or an expression of the form $P(t_1, ..., t_n)$ where $P \in \mL$ is any $n$-ary predicate symbol and $t_1, ..., t_n$ are terms of $\mL$ over $\mC$. The \la{formulas } $\mF(\mL, \mC)$ of $\mL$ over $\mC$ form the smallest set of expressions containing the atomic formulas and closed under the  formation rule: if $A, B$ are formulas so are the expressions $(A \Rightarrow B)$ and $(\forall A)$. We shall follow the usual conventions to eliminate parentheses.

In the following we assume $P, Q, R \in \mL$ are predicate symbols, $f, g, h \in \mL$ are function symbols,  $x, y, z \in X = \{x_1, x_2, ...\}$ are variables, $c, c_1, c_2, ... \in \mC$ are parameters, $s, s_1, s_2, ..., t, t_1, t_2, ... \in \mT(\mL, \mC)$ are terms over $\mC$, and $A, B, C \in \mF(\mL, \mC)$ are formulas over $\mC$.

Define a mapping $\mT(\mL, \mC) \times \mT(\mL, \mC)^{\omega} \ar \mT(\mL, \mC)$ inductively as follows:
\\ (i) $x_i[t_1, t_2, ...] = t_i$.
\\ (ii) $c[t_1, t_2, ...] = c$.
\\ (iii) $f(t_1, ..., t_n)[s_1, s_2, ....] = f(t_1[s_1, s_2, ...], t_2[s_1, s_2, ...], ...])$.

Define a mapping $\tau: \mF(\mL, \mC) \times \mT(\mL, \mC)^{\omega} \ar \mF(\mL, \mC)$ inductively as
follows:
\\ 1. $P(s_1, ..., s_n)[t_1, t_2, ...] = P(s_1[t_1, t_2, ...], ...,
s_n[t_1, t_2, ...])$.
\\ 2. $(A \Rightarrow B)[t_1, t_2, ...] = (A[t_1, t_2, ...] \Rightarrow
B[t_1, t_2, ...])$.
\\ 3. $(\forall  A)[t_1, t_2, ...] = \forall (A[x_1, t_1^+, t_2^+, ....])$ where $t_i^+ = t_i[x_2, x_3, ...]$.

Let $\mT(\mC) := \mT(\mL, \emptyset)$ and $\mF(\mL) := \mF(\mL, \emptyset)$. Note that $\mT(\mC) \subset \mT(\mL, \mC)$ and $\mF(\mL) \subset \mF(\mL, \mC)$ for any  $\mC$. Define $\neg A := (A \Rightarrow \FF)$, $x \approx y := \approx(x, y) $ and $A[t/x_i] := A[x_1, ..., x_{i-1}, t, x_{i+1}, ...]$. Define $\forall^0 A = A$ and $\forall^n A = \forall (\forall^{n-1} A)$ for $n > 1$ inductively.

\ps{(cf. \cite{luo:3}) 1. $\mT(\mL, \mC)$ is a locally finite clone, which is a free algebra over the free basis $X \cup \mC$ with function symbols as the  signature.
\\ 2. $\mF(\mL, \mC)$ is a locally finite free predicate algebra over the clone $\mT(\mL, \mC)$ generated by atomic formulas.

}

Suppose $D$ is a term or formula. We say $D$ is \la{independent} of a variable $x_i$ if $D = D[x_{i+1}/x_i]$. If $D$ is not independent of $x$ then we say that $x$ is \la{free} in $D$. The set of free variables in $D$ is always finite.  Let $D^+ := D[x_2, x_3, ...]$ and $D^- := D[x_1, x_1, x_2, ...]$. Then $(D^+)^- = D$. We say $D$ has a \la{rank} $n \ge 0$ if $D = D[x_1, ..., x_{n-1}, x_n, x_n, ...]$.

 Denote by $\mF_n(\mL, \mC)$ the set of formulas with a rank $n \ge 0$.  If $A, B \in \mF_n(\mL, \mC)$ then $ (A \Rightarrow B), (\forall A), (\forall x_i)A \in \mF_n(\mL, \mC)$. If $n > 0$ and $A \in \mF_n(\mL, \mC)$ then $(\forall A), (\forall x_n)A \in \mF_{n-1}(\mL, \mC)$. A \la{sentence} is a formula with a rank $0$.  If a formula $A$ has a rank $n > 0$ then $\forall^n A = (\forall x_n)...(\forall x_1)A$ is a sentence.

For any variable $x_i$ let $(\forall x_i)A := \forall(A[x_2, x_3, ..., x_i, x_1, x_{i+2}, ...])$.

\lm{1. $(\forall x)A$ is independent of $x$.
\\ 2. $(\forall x)A = (\forall y)(A[y/x])$ if $A$ is independent of $y$.
\\ 3. $((\forall x_i)A)[t_1, t_2, ...] = (\forall y)(A[t_1, ..., t_{i-1}, y, t_{i+1}, ...])$ if $t_j$ is independent of $y$ for any $j \ne i$ such that $x_j$ is free in $A$.}

\lm{
1. $\forall A = (\forall x_i)(A[x_i, x_1, x_2, ...])$ if $A$ is independent of $x_{i+1}$.
\\ 2. $(\forall x)A = \forall (A^+) $ if $A$ is independent of $x$.
\\ 3. $(\forall x_1)A = (\forall A)^+$ and $\forall A = ((\forall x_1)A)^-$.
\\ 4. $(\forall x_i)A = ((\forall x_1)(A[x_2, x_3, ..., x_i, x_1, x_{i+2}, ...]))^-$.

}

 A \la{perfect  valuation} (or \la{Henkin  valuation}) of $\mL$ (over $\mC$) is a subset  $U$ of $\mF(\mL, \mC)$ satisfying the following conditions for any $A, B \in \mF(\mL, \mC)$ and $t, t_1, t_2, ... \in \mT(\mL, \mC)$:
\\ 1. $\FF \notin U$.
\\ 2. $(A \Rightarrow B) \in U$ iff $A \notin U$ or $B \in U$.
\\ 3. $(\forall A) \in U$ iff $A[t, x_1, x_2, ...] \in U$ for any term $t$. (or equivalently, for any $x$, $(\forall x) A \in U$ iff $A[t/x] \in U$  for any term $t$).
\\ If $\mL$ is a language with equality then we also assume that the following conditions are satisfied.
\\ 4. $\forall^n(x \approx x) \in U$ for any $n \ge 0$.
\\ 5. $\forall^n(x \approx y \Rightarrow (A \Rightarrow A[y/x])) \in U$ for any $n \ge 0$.

 Denote by $\mA(\mL, \mC)$ the set of atomic formulas of $\mL$ over $\mC$. A subset $E$ of $\mA(\mL, \mC)$ such that $\FF \notin E$ is called an \la{atomic valuation} of $\mL$ over $\mC$. Since $\mF(\mL, \mC)$ is generated by atomic formulas under operations $\Rightarrow$ and $\forall$ inductively, a perfect valuation $U$ of $\mL$ is uniquely determined by the atomic valuation $U \cap \mA(\mL, \mC)$. Conversely any atomic valuation determines a perfect valuation for any language $\mL$ without equality.

 If $U \subset \mF(\mL, \mC)$ is a perfect valuation and $t_1, t_2, ... \in \mT(\mL, \mC)$, let $U_{(t_1, t_2, ...)} := \{A \in \mF(\mL, \mC) \ | \ A[t_1, t_2, ...] \in U\}$. A subset $W$ of $\mF(\mL, \mC)$ is called a \la{valuation} of $\mL$ (over $\mC$) if $W = U_{(t_1, t_2, ...)}$ for some $U$ and $t_1, t_2, ... $ as above.  A subset $V$ of $\mF(\mL)$ is called a \la{logical valuation} of $\mL$ if there is a valuation $W$ of $\mL$ (over some set $\mC$ of parameters) such that $V = \mF(\mL) \cap W$. A subset of $\mF(\mL)$ is called   \la{logically closed} if it is an intersection of logical valuations of $\mL$. A formula $A \in \mF(\mL)$ is called \la{logically valid} if it belongs to any logical valuation of $\mL$.

A \la{structure} for $\mL$ is a pair $\mM = (M, \gamma)$ where $M$ is a set and $\gamma$ is an operation with domain $\mL$ such that
\\ (i) if $f \in \mL$ is an $n$-ary function symbol then, then $\gamma(f): M^n \ar M$.
\\ (ii)  if $P \in \mL$ is an $n$-ary predicate symbol then, then $\gamma(P) \subset M^n$.
\\ (iii) $\gamma(\FF) = \emptyset$.
\\ (iv) $(m_1, m_2) \in \gamma(\approx)$ iff $m_1 = m_2$.

 Any structure $\mM = (M, \gamma)$ determines a left algebra $M$ over the clone $\mT(\mL, M)$ such that $f(x_1, ..., x_n)[m_1, m_2, ...] = \gamma(f)(m_1, ..., m_n)$ and $m[m_1, m_2, ...] = m$ for any elements $m, m_1, m_2, ... \in M$. Any fixed sequence $m_1,  m_2, ... \in M$ then determines a perfect valuation $U$ of $\mL$ over $M$ such that $P(t_1, ..., t_n) \in U$ iff $(t_1[m_1, m_2, ...], ..., t_n[m_1, m_2, ...]) \in \gamma(P)$. Conversely, any perfect valuation $U$ of $\mL$ over $\mC$ determines a structure $(\mT(\mL, \mC), \gamma)$ such that $\gamma(f)(t_1, ..., t_n) = f(t_1, ..., t_n)$ and $(t_1, ..., t_n) \in \gamma(P)$ iff $P(t_1, ..., t_n) \in U$.

Suppose $A, B, C \in \mF(\mL)$. The following formulas are called \la{(first order) axioms}:
\\ A1. $A \Rightarrow (B \Rightarrow A)$.
\\ A2. $(A \Rightarrow (B \Rightarrow C)) \Rightarrow ((A \Rightarrow
B) \Rightarrow (A \Rightarrow C))$.
\\ A3. $\neg (\neg A) \Rightarrow A$.
\\ A4. $\forall (A \Rightarrow B) \Rightarrow (\forall A \Rightarrow \forall B)$.
\\ A5. $\forall A \Rightarrow A[t, x_1, x_2, ...]$ for any term $t \in \mT(\mL)$.
\\ A6.  $A \Rightarrow \forall (A^+)$.
\\ If $\mL$ is a language with equality then the following formulas are also axioms:
\\ A7. $x \approx x$.
\\ A8. $x \approx y \Rightarrow (A \Rightarrow A[y/x])$.
\\ Furthermore if $A$ is an axiom then $\forall A$ is an axiom.

Note that A4-A6 are equivalent to the following A4$'$-A6$'$ for any variable $x$:
\\ A4$'$. $(\forall x) (A \Rightarrow B) \Rightarrow ((\forall x) A
\Rightarrow (\forall x) B)$.
\\ A5$'$. $(\forall x)A \Rightarrow A[t/x]$ for any term $t \in \mT(\mL)$.
\\ A6$'$. $A \Rightarrow (\forall x)A$ if $A$ is independent of $x$.

A \la{(first order) filter} of $\mF(\mL)$ is a subset $I$ of $\mF(\mL)$ containing all the axioms such that if $A, (A \Rightarrow B) \in I$ then $B \in I$.
If $S$ is any subset of $\mF(\mL)$ denote by $Con(S)$ the intersection of all logical valuations containing $S$. Denote by $Ded(S)$ the intersection of all filters containing $S$. Write $S \vDash A$ if $A \in Con(S)$, $S \vdash A$ if $A \in Ded(S)$,  and $S \nvdash A$ if $A \notin Ded(S)$.

\te{1. (Soundness Theorem) Any logically closed set is a filter.
\\ 2. (Completeness Theorem) Any filter is logically closed.
\\ 3. $S \vDash A$ iff $S \vdash A$ (i.e. $Con(S) = Ded(S)$). }

A \la{theory} of $\mL$ is a set $T$ of sentences. We say a theory $T$ is \la{consistent} if there is no formula $A$ such that $T \vdash A$ and $T \vdash \neg A$. A theory $T$ is \la{complete} if for any sentence $A$, we have $T \vdash A$ iff $T \nvdash \neg A$.

Let $\mL_a$ be a language with equality and function symbols $\bf 0$, $'$,  $\cdot$ and $+$.

Let $T_a$ be the theory consisting of the following sentences of $\mL_a$:
\\ (S1) $(\forall x) (\neg ({\bf 0}  \approx x'))$.
\\ (S2) $(\forall x)(\forall y)((x' \approx y') \Rightarrow (x \approx y))$.
\\ (S3) $ (\forall x)(x + {\bf 0} \approx  x)$.
\\ (S4) $ (\forall x)(\forall y)(x + y' \approx (x + y)')$.
\\ (S5) $ (\forall x)(x \cdot \bf 0 \approx \bf 0)$.
\\ (S6) $ (\forall x)(\forall y)(x \cdot y' \approx x \cdot y + x)$.
\\ (S7) $(\forall x_n)...(\forall x_1) (A[{\bf 0}/x] \Rightarrow ( (\forall x) (A \Rightarrow A[x'/x]) \Rightarrow  (\forall x)A))$ for any formula $A$ of rank $n > 0$.

\te{(Incompleteness Theorem). Assume $T_a$ is consistent. Then it is not complete.
}

\end{document}